\documentclass[12pt]{amsart}
\usepackage{amssymb}
\usepackage{isolatin1}    
\usepackage{kuvio}

\newtheorem{thm}{Theorem}[section]
\newtheorem{lem}[thm]{Lemma}

\newtheorem{prop}[thm]{Proposition}
\newtheorem{ex}[thm]{Example}

\newtheorem*{prob*}{Open problem}

\theoremstyle{definition}

\newtheorem{defi}[thm]{Definition}

\theoremstyle{remark}

\newtheorem{rem}[thm]{Remark}
\newtheorem*{rem*}{Remark}

\DeclareMathOperator{\s}{span}

\newcommand{\kringel}{\mathbin{\raise1pt\hbox{$\scriptstyle\circ$}}} 
\newcommand{\pkt}{\mathbin{\raise0pt\hbox{$\scriptstyle\bullet$}}}

\newcommand{\C}{\mathbb{C}}

\newcommand{\N}{\mathbb{N}}
\renewcommand{\P}{\mathbb{P}}

\newcommand{\tr}{\mathop{\rm tr}}
\newcommand{\ad}{\mathop{\rm ad}}

\newcommand{\Der}{\mathop{\rm Der}}
\newcommand{\diag}{\mathop{\rm diag}}

\newcommand{\La}{\mathfrak{a}}

\newcommand{\Lg}{\mathfrak{g}}
\newcommand{\Lh}{\mathfrak{h}}

\newcommand{\Ll}{\mathfrak{l}}
\newcommand{\Ln}{\mathfrak{n}}

\newcommand{\Lr}{\mathfrak{r}}
\newcommand{\Ls}{\mathfrak{s}}
\newcommand{\Lt}{\mathfrak{t}}

\newcommand{\CC}{\mathcal{C}}
\newcommand{\CL}{\mathcal{L}}
\newcommand{\CN}{\mathcal{N}}

\newcommand{\CR}{\mathcal{R}}

\newcommand{\abs}[1]{\lvert#1\rvert}

\newcommand{\al}{\alpha}
\newcommand{\be}{\beta}

\newcommand{\la}{\lambda}

\newcommand{\ov}{\overline}

\newcommand{\ra}{\rightarrow}  
\newcommand{\ch}{\checkmark} 

\renewcommand{\phi}{\varphi}

\begin{document}


\title[Degenerations of nilpotent Lie algebras]{Degenerations of
$\bf{7}$-dimensional nilpotent Lie algebras.} 

\author[D. Burde]{Dietrich Burde}
\address{Fakult\"at f\"ur Mathematik\\
Universit\"at Wien\\
  Nordbergstr. 15\\
  1090 Wien \\
  Austria}
\date{\today}
\email{dietrich.burde@univie.ac.at}

\begin{abstract}
We study the varieties of Lie algebra laws and their
subvarieties of nilpotent Lie algebra laws. 
We classify all degenerations of (almost all) five-step and six-step nilpotent 
seven-dimensional complex Lie algebras. 
One of the main tools is the use of trivial and adjoint cohomology of
these algebras. In addition, we give some new results on the varieties 
of complex Lie algebra laws in low dimension.
\end{abstract}

\maketitle

\section{Introduction} 

Let $\Lg$ be an $n$-dimensional vector space over a field $k$
and consider the set $\CL_n(k)$ of all possible Lie brackets $\mu$
on $\Lg$. This is an algebraic subset of the variety   
$\Lambda^2 \Lg^{\ast} \otimes \Lg$ of all alternating
bilinear maps from $\Lg\times \Lg$ to $\Lg$. Indeed, for a fixed
basis $(x_1,\ldots,x_n)$ of $\Lg$ the Lie bracket $\mu$ is determined
by the point $(c_{ijr})\in k^{n^3}$ of structure constants with
\begin{equation*}
\mu(x_i,x_j) = \sum_{r=1}^n c_{ijr} x_r
\end{equation*}
satisfying the polynomial conditions
\begin{align*}\label{jac}
c_{ijr} + c_{jir} & =0, \\
\sum_{r=1}^n (c_{ijr} c_{lrs}+c_{jkr} c_{irs}+c_{kir} c_{jrs}) & = 0,
\quad 1\le i<j<k\le n,\; 1\le s\le n
\end{align*}
The variety $\CL_n(k)$ is often called the {\it variety of Lie algebra laws}.
The general linear group $GL_n(k)$ acts on $\CL_n(k)$ by base change:
\begin{equation}\label{action}
(g\cdot\mu)(x,y)=g(\mu(g^{-1}x, g^{-1}y)),\quad g\in GL_n(k),\; x,y\in \Lg
\end{equation}
One denotes by $O(\mu)$ the orbit of $\mu$ under the action of $GL_n(k)$, and
by $\ov{O(\mu)}$ the closure of the orbit with respect to the Zariski
topology.
The orbits in $\CL_n(k)$ correspond to isomorphism classes of $n$-dimensional
Lie algebras. However, the orbit space is no longer an algebraic set. It makes
sense to take out the zero point and to view 
$$(\CL_n(k)\setminus {0})/GL_n(k)=\P(\CL_n(k))/PGL_n(k)$$ 
as the {\it moduli space}.\\
There are many questions on the structure of the varieties $\CL_n(k)$. In
particular one is interested in the irreducible components of $\CL_n(k)$ and
in the open orbits. A Lie algebra law $\mu\in \CL_n(k)$ is called {\it rigid},
if its orbit $O(\mu)$ is open in $\CL_n(k)$. In that case the corresponding Lie algebra
$\Lg$ is algebraic and does not admit any non-trivial deformation \cite{CA2}.
On the other hand $H^2(\Lg,\Lg)=0$ implies that $\mu$ is rigid. The converse
does not hold in general. The following result (see \cite{CA1})
gives the number of components and open orbits in $\CL_n(\C)$ in low dimensions:

\begin{prop}
Let $r(n)$ denote the number of irreducible components in $\CL_n(\C)$ 
and $s(n)$ the number of open orbits. Then it holds 
$(r(1),\ldots,r(7))=(1,1,2,4,7,17,49)$ and
$(s(1),\ldots,s(7))=(1,1,1,2,3,6,14)$.
\end{prop}
\vspace{0.5cm}
These numbers grow very fast in $n$. One has the following estimates
for $n$ big enough \cite{KN}:
\begin{equation*}
e^{n/4} < s(n) < r(n) < 2^{n^4/6}
\end{equation*}   

In studying the orbit closures the concept of Lie algebra degenerations
is of great interest.

\begin{defi}
We say that $\mu$ is a degeneration of $\la$ in $\CL_n(k)$ if
$\mu \in \ov{O(\la)}$. In that case we also say that $\la$ degenerates
to $\mu$, which is denoted by $\la \ra_{\rm deg} \mu$.
\end{defi}

Let $\CC$ be an irreducible component of $\CL_n(k)$ containing $\mu$.
Then also $O(\mu)\subset \CC$. Since $\CC$ is closed relative to the
Zariski topology, the orbit closure $\ov{O(\mu)}$ is contained in
 $\CC$. Hence any irreducible component containing $\mu$ also contains
all degenerations of $\mu$.

\begin{prop}
Degeneration defines an order relation on the orbit space of $n$-dimensional
Lie algebra laws by $O(\mu)\le O(\la) \iff \mu \in \ov{O(\la)}$
\end{prop}

\begin{proof}
The relation is clearly reflexive. The transitivity follows from the fact that
$O(\la)\subseteq \ov{O(\mu)} \iff \ov{O(\la)}\subseteq \ov{O(\mu)}$.
Finally, antisymmetry follows from the fact, that any orbit in this case
is open in its closure.
\end{proof}

A degeneration is called {\it trivial} if $\la \cong \mu$, that is, if
 $\mu \in O(\la)$. Note that $\la \ra_{\rm deg} \mu$ and 
$\mu \ra_{\rm deg} \nu$ imply that $\la \ra_{\rm deg} \nu$. That is the
transitivity of the above order relation.

\begin{rem}
The concept of degenerations was first introduced by theoretical physicists
in the special case of {\it contractions} \cite{IW}.
Often the limit procedures considered in physics can be described by Lie algebra
contractions. As an example, classical mechanics is a limit of quantum mechanics
given by the contraction $\Lh \ra_{\rm deg} \Lt_{2n+1}$, where $\Lh$ is the
Weyl-Heisenberg algebra and $\Lt_{2n+1}$ is the abelian Lie algebra of the
same dimension.
\end{rem}

It is known that over the real or complex numbers the Zariski closure of
an orbit coincides with the orbit closure relative to the usual metric
topology. The definition of Lie algebra degeneration can be refrased
so that the relation to Lie algebra deformations can be made apparent \cite{GRH}:

\begin{prop}
Let $k$ be an algebraically closed field and $\Lg$ and $\Lh$ two
$n$-dimensional Lie algebras over $k$. Then $\Lh$ is a 
degeneration of $\Lg$ if and only there exists a discrete
valuation algebra $A$ over $k$ such that its field of fractions $K$
is a function field of dimension $1$, and if there is a Lie algebra
$\La$ over $A$ of dimension $n$ such that
\begin{align*}
\La \otimes_A K & \cong \Lg \otimes_k K\\
\La \otimes_A k & = \Lh 
\end{align*} 
\end{prop} 

Often a degeneration can be realized by a one-parameter subgroup
$\{g_t\}$ of $GL_n(k)$, see \cite{BU9}
\begin{defi}
A degeneration $\la \ra_{\rm deg} \mu$ is called a  
{\it one-parameter subgroup degeneration}, or 1-PSG, if it can be
realized by a group homomorphism $g:k^* \ra GL_n(k), \; t \mapsto g_t$ such that
$\mu\cong\lim_{t \to 0} g_t\cdot\la$.                        
\end{defi}

The notion of a 1-PSG degeneration does not depend on the choice of a basis.

\begin{ex}
Let $\la_0\in \CL_n(k)$ be the law corresponding to the abelian Lie algebra,
i.e., $\la_0(x,y)=0$, and $g_t=t^{-1}I_n$.
Then $\la \ra_{\rm deg} \la_0$ for all $\la \in \CL_n(k)$:
$$(g_t \cdot \la)(x,y)=t^{-1}\la(tx,ty)=t\la(x,y)$$
Indeed, the limit of $g_t \cdot \la$ for $t\to 0$ equals 
$\la_0$. Hence every Lie algebra degenerates to the abelian Lie algebra
of the same dimension by a 1-PSG degeneration. For some Lie algebras, such as
$\Lh_3\oplus k^m$, where $\Lh_3$ is the $3$-dimensional Heisenberg Lie algebra,
this is the only possible degeneration, see \cite{LA}.
\end{ex}

Given two Lie algebra laws $\la,\mu \in \CL_n(k)$ it is sometimes quite
difficult to see whether there exists a degeneration $\la \ra_{\rm deg} \mu$.
It is helpful to obtain some necessary conditions for the existence of a 
degeneration.
In some sense one can say that $\la \ra_{\rm deg} \mu$ implies that
$\mu$ is ``more abelian'' than $\la$. A much finer condition is that
the dimensions of the cohomology spaces cannot decrease. 

\begin{prop}\label{inv}
Let $\la \ra_{\rm deg} \mu$ a non-trivial degeneration. Then we have
for all $j\in \N_0$:\\
\begin{align*}
\dim O(\la) & > \dim O(\mu)\\
\dim \Der \la  & < \dim \Der \mu\\
\dim [\mu,\mu] & \le \dim [\la,\la]\\
\dim Z(\la) & \le \dim Z(\mu) \\
\dim Z^j(\la) & \le \dim Z^j(\mu) \\
\dim Z^j(\la,\la) & \le \dim Z^j(\mu,\mu)\\
\dim H^j(\la) & \le \dim H^j(\mu) \\
\dim H^j(\la,\la) & \le \dim H^j(\mu,\mu)\\
\end{align*}
\end{prop}

\begin{proof}
These inequalities are well known. I have not seen the ones
on cohomology in the literature yet. So let us repeat the argument. 
It is clear that we have
$\dim Z^j(\la,\la) \le \dim Z^j(\mu,\mu)$ for $j\in \N_0$.
Let $d: C^j(\la,\la)\ra C^{j+1}(\la,\la)$ be the coboundary operator
of the standard complex for the Lie algebra cohomology.
Using the dimension formula for the linear map $d$ we have
$$\dim H^j(\la,\la)=\dim Z^j(\la,\la)-\dim C^{j+1}(\la,\la)+
\dim Z^{j+1}(\la,\la)$$
But that implies $\dim H^j(\la,\la) \le \dim H^j(\mu,\mu)$, since
$\dim C^{j+1}(\la,\la)=\dim C^{j+1}(\mu,\mu)$.
The same argument applies for the cohomology with trivial coefficients.
\end{proof}

\section{Degenerations in dimension $6$}

Denote by $\CN_n(k)$ the subvariety of $\CL_n(k)$ consisting
of $n$-dimensional nilpotent Lie algebra laws. It is known that
the varieties $\CN_n(\C)$ are irreducible for $n\le 6$ and
reducible for all $n\ge 11$ \cite{AGG} and $n=7,8,9$. For $n\le 6$ 
all degenerations in $\CN_n(\C)$ are known, see \cite{GRH},\cite{SEE}.
We will shortly summarize the results. \\[0.5cm] 
Let $\Ln_3(\C)$ denote the $3$-dimensional Heisenberg Lie algebra. We have
$$ 
\CN_3(\C)=\ov{O(\Ln_3(\C))}=O(\Ln_3(\C))\cup O(\C^3)
$$
the only (non-trivial) degeneration being $\Ln_3(\C)\ra_{\rm deg} \C^3$.
For $n=4$ we have
$$ 
\CN_4(\C)=\ov{O(\Ln_4(\C))}=O(\Ln_4(\C))\cup O(\Ln_3(\C)\oplus \C)\cup O(\C^3)
$$
where $\Ln_4(\C)$ is the standard graded filiform Lie algebra of dimension $4$.
The degenerations are given by

$$
\Dg
\Ln_4 & \rTo & \Ln_3\oplus \C & \rTo & \C^4 \\
\endDg
$$

For $n=5$ we have the following classification of all orbits
in $\CN_5(\C)$:
\vspace*{0.5cm}
\begin{center}
\begin{tabular}{c|c}
 $\Lg$ & Lie brackets \\
\hline
$\C^5$ & $-$\\
$\Ln_3(\C)\oplus \C^2$ & $[e_1,e_2]=e_3$\\
$\Ln_4(\C)\oplus \C$ & $[e_1,e_2]=e_3, [e_1,e_3]=e_4$\\
$\Lg_{5,6}(\C)$ & $[e_1,e_2]=e_3, [e_1,e_3]=e_4, [e_1,e_4]=e_5, [e_2,e_3]=e_5$ \\
$\Lg_{5,5}(\C)$ & $[e_1,e_2]=e_3, [e_1,e_3]=e_4, [e_1,e_4]=e_5$ \\
$\Lg_{5,4}(\C)$ & $[e_1,e_2]=e_3, [e_1,e_3]=e_4, [e_2,e_3]=e_5$ \\
$\Lg_{5,3}(\C)$ & $[e_1,e_2]=e_4, [e_1,e_4]=e_5, [e_2,e_3]=e_5$ \\
$\Lg_{5,2}(\C)$ & $[e_1,e_2]=e_4, [e_1,e_3]=e_5$ \\
$\Lg_{5,1}(\C)$ & $[e_1,e_3]=e_5, [e_2,e_4]=e_5$ \\
\end{tabular}
\end{center}   
\vspace*{0.5cm}
The degenerations in $\CN_5(\C)$ have been classified in \cite{GRH}. The 
Hasse diagram is given by:
 $$
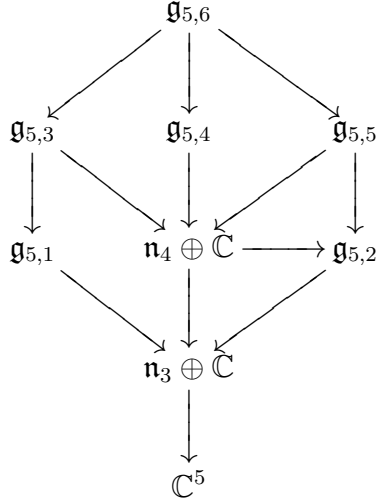
\Dg
             &      & \Lg_{5,6} \aTo (-2,-2) \aTo (2,-2) &       &      \\
  &      &  \dTo     &       &   \\
\Lg_{5,3} \aTo (2,-2)   &      & \Lg_{5,4} &       & \Lg_{5,5}  \aTo (-2,-2) \\
   \dTo     &      & \dTo              &      &  \dTo \\
\Lg_{5,1}\aTo (2,-2) & \phantom{---} & \Ln_4\oplus \C & \rTo &  \Lg_{5,2} \aTo (-2,-2)\\
             &      &  \dTo             &      &  \\
             &      &  \Ln_3\oplus\C    &      &  \\
             &      &   \dTo            &      &  \\
             &      &   \C^5            &      &  \\
\endDg
$$   
The Lie algebra $\Lg_{5,6}(\C)$ is on top of the diagram. It is rigid
in $\CN_5(\C)$, hence 
$$\CN_5(\C)=\ov{O(\Lg_{5,6}(\C))}$$
For $n=6$ the degeneration diagram becomes very complicated, see \cite{SEE}.
Restricting ourselfs to filiform Lie algebras of dimension six the picture 
becomes much easier. The classification is given by
\vspace*{0.5cm}
\begin{center}
\begin{tabular}{c|c}
 $\Lg$ & Lie brackets \\
\hline
$\Lg_{6,A}$ & $[e_1,e_i]=e_{i+1},\, 2\le i\le 5 $ \\
$\Lg_{6,B}$ & $[e_1,e_i]=e_{i+1},\; [e_2,e_3]=e_6$ \\
$\Lg_{6,C}$ & $[e_1,e_i]=e_{i+1},\; [e_2,e_5]=e_6, [e_3,e_4]=-e_6$ \\
$\Lg_{6,D}$ & $[e_1,e_i]=e_{i+1},\; [e_2,e_3]=e_5, [e_2,e_4]=e_6$ \\
$\Lg_{6,E}$ & $[e_1,e_i]=e_{i+1},\; [e_2,e_3]=e_5, [e_2,e_4]=e_6,
 [e_2,e_5]=e_6,[e_3,e_4]=-e_6$ \\
\end{tabular}
\end{center}   
\vspace*{0.5cm}
The degenerations among filiform Lie algebras in $\CN_6(\C)$ are given by:
 $$
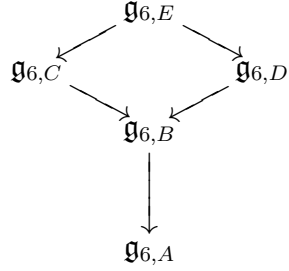
\Dg
                  & \Lg_{6,E} \aTo (-1,-1) \aTo (1,-1)   &      \\
\Lg_{6,C} \aTo (1,-1)    &    &     \Lg_{6,D}  \aTo (-1,-1) \\
                         & \Lg_{6,B}  &   \\
                         & \dTo & \\
              &  \Lg_{6,A}    &     \\
\endDg
$$   
The Lie algebra $\Lg_E$ is rigid in $\CN_6(\C)$, hence 
$\CN_6(\C)=\ov{O(\Lg_{6,E})}$.

\section{Degenerations in dimension 7}

The classification of all degenerations of complex nilpotent
Lie algebras of dimension $7$ is an enormous task. We cannot consider
all Lie algebras here. Instead we restrict ourselfs to the subset of 
indecomposable Lie algebras of nilpotency
class $5$ and $6$. From the list given in \cite{MAG} we conclude
that these are the following Lie algebras. We add the notation used
in \cite{MAG}. 

\begin{align*}
\Lg_{I}(\al) & = \Lg_{7,1.1(i_{\la})}, \; \al\ne 0 \hspace{13cm} \\
[x_1,x_i] & = x_{i+1}; \; 2\le i\le 6;\quad
[x_2,x_3] =x_5; \; [x_2,x_4]=x_6; \; [x_2,x_5]=(1-\al)x_7;\; [x_3,x_4]=\al x_7.
\end{align*}
\begin{align*}
\Lg_{F} \; & =  \Lg_{7,0.1} \hspace{16cm} \\
[x_1,x_i] & = x_{i+1}; \; 2\le i\le 6 ;\quad
[x_2,x_3] =x_6; \; [x_2,x_4]=x_7; \; [x_2,x_5]=x_7;\; [x_3,x_4]=-x_7.
\end{align*}
\begin{align*}
\Lg_{H} \; & =  \Lg_{7,0.2} \hspace{16cm} \\
[x_1,x_i] & = x_{i+1}; \; 2\le i\le 6 ;\quad 
[x_2,x_3] =x_5 + x_7; \; [x_2,x_4]=x_6; \; [x_2,x_5]=x_7.
\end{align*}
\begin{align*}
\Lg_{1}(\la) & = \Lg_{7,0.4(\la)} \hspace{16cm} \\
[x_1,x_2] & = x_3; \;  [x_1,x_3] = x_4; \; [x_1,x_4] = x_6+\la x_7; \; 
[x_1,x_5] = x_7; \; [x_1,x_6] = x_7; \; \\
[x_2,x_3] & =x_5; \; [x_2,x_4]=x_7; \; [x_2,x_5]=x_6;\; [x_3,x_5]=x_7.
\end{align*}
\begin{align*}
\Lg_{2} \; & = \Lg_{7,0.5} \hspace{16cm} \\
[x_1,x_2] & = x_3; \; [x_1,x_3] = x_4; \; [x_1,x_4] = x_6+ x_7; \; 
[x_1,x_6] = x_7; \; \\
[x_2,x_3] & =x_5; \; [x_2,x_5]=x_6;\; [x_3,x_5]=x_7.
\end{align*}
\begin{align*}
\Lg_{3} \; & = \Lg_{7,0.6} \hspace{16cm} \\
[x_1,x_2] & = x_3; \; [x_1,x_3] = x_4; \; [x_1,x_4] = x_7; \; 
[x_1,x_5] = x_6; \;[x_1,x_6] = x_7; \; \\
[x_2,x_3] & =x_5; \; [x_2,x_4]=x_6;\;  [x_2,x_5]=x_7;\; [x_3,x_4]=x_7.
\end{align*}
\begin{align*}
\Lg_{4} \; & = \Lg_{7,0.7} \hspace{16cm} \\
[x_1,x_2] & = x_3; \; [x_1,x_3] = x_4; \; [x_1,x_4] = x_5; \; 
[x_1,x_6] = x_7; \; \\
[x_2,x_3] & =x_5+x_6; \; [x_2,x_4]=x_7;\;  [x_2,x_5]=x_7;\; [x_3,x_4]=-x_7.
\end{align*}
\begin{align*}
\Lg_{5} \; & = \Lg_{7,0.8} \hspace{16cm} \\
[x_1,x_2] & = x_4; \; [x_1,x_3] = x_7; \; [x_1,x_4] = x_5; \; 
[x_1,x_5] = x_6; \; \\
[x_2,x_3] & =x_6; \; [x_2,x_4]=x_6;\;  [x_2,x_6]=x_7;\; [x_4,x_5]=-x_7.
\end{align*}
\begin{align*}
\Lg_{6} \; & = \Lg_{7,1.1 (iii)} \hspace{16cm} \\
[x_1,x_2] & = x_3; \; [x_1,x_3] = x_4; \; [x_1,x_4] = x_5; \; 
[x_1,x_5] = x_6; \; \\
[x_2,x_3] & =x_5; \; [x_2,x_4]=x_6;\;  [x_2,x_5]=-x_7;\; [x_3,x_4]=x_7.
\end{align*}
\begin{align*}
\Lg_{7} \; & = \Lg_{7,1.1 (v)} \hspace{16cm} \\
[x_1,x_3] & = x_4; \; [x_1,x_4] = x_5; \; [x_1,x_5] = x_6; \; 
[x_1,x_6] = x_7; \; \\
[x_2,x_3] & =x_5; \; [x_2,x_4]=x_6;\;  [x_2,x_5]=x_7;\; [x_3,x_4]=-x_7.
\end{align*}
\begin{align*}
\Lg_{C} \; & = \Lg_{7,1.1 (ii)} \hspace{16cm} \\
[x_1,x_i] & = x_{i+1}; \; 2\le i\le 6;\quad [x_2,x_5] =x_7; \; [x_3,x_4]=-x_7.
\end{align*}
\begin{align*}
\Lg_{G} \; & = \Lg_{7,1.1 (i_0)} \hspace{16cm} \\
[x_1,x_i] & = x_{i+1}; \; 2\le i\le 6 ;\quad 
[x_2,x_3] =x_5; \; [x_2,x_4] =x_6; \; [x_2,x_5] =x_7.
\end{align*}
\begin{align*}
\Lg_{E} \; & = \Lg_{7,0.3} \hspace{16cm} \\
[x_1,x_i] & = x_{i+1}; \; 2\le i\le 6 ;\quad 
[x_2,x_3] =x_6+x_7; \; [x_2,x_4] =x_7.
\end{align*}
\begin{align*}
\Lg_{8} \; & = \Lg_{7,1.01(i)} \hspace{16cm} \\
[x_1,x_3] & = x_4; \; [x_1,x_4] = x_5; \; [x_1,x_5] = x_6; \; 
[x_1,x_6] = x_7; \; \\
[x_2,x_3] & =x_5+x_7; \; [x_2,x_4]=x_6;\;  [x_2,x_5]=x_7.
\end{align*}
\begin{align*}
\Lg_{9} \; & = \Lg_{7,1.02} \hspace{16cm} \\
[x_1,x_2] & = x_3; \; [x_1,x_3] = x_4+x_6; \; [x_1,x_5] = x_6; \; 
[x_1,x_6] = x_7; \; \\
[x_2,x_3] & =x_5; \; [x_2,x_4]=x_6;\;  [x_3,x_4]=x_7.
\end{align*}
\begin{align*}
\Lg_{10} \; & = \Lg_{7,1.03} \hspace{16cm} \\
[x_1,x_2] & = x_3; \; [x_1,x_3] = x_4; \; [x_1,x_4] = x_5; \; 
[x_1,x_6] = x_7; \; \\
[x_2,x_3] & =x_6; \; [x_2,x_4]=x_7;\;  [x_2,x_5]=x_7;\; [x_3,x_4]=-x_7.
\end{align*}
\begin{align*}
\Lg_{11} \; & = \Lg_{7,1.1(iv)} \hspace{16cm} \\
[x_1,x_2] & = x_3; \; [x_1,x_3] = x_4; \; [x_1,x_5] = x_6; \; 
[x_1,x_6] = x_7; \; \\
[x_2,x_3] & =x_5; \; [x_2,x_4]=x_6;\;  [x_2,x_5]=x_7;\; [x_3,x_4]=x_7.
\end{align*}
\begin{align*}
\Lg_{12} \; & = \Lg_{7,1.1(vi)} \hspace{16cm} \\
[x_1,x_2] & = x_3; \; [x_1,x_3] = x_4; \; [x_1,x_4] = x_5; \; 
[x_1,x_6] = x_7; \; \\
[x_2,x_3] & =x_5; \; [x_2,x_5]=x_7;\;  [x_3,x_4]=-x_7.
\end{align*}
\begin{align*}
\Lg_{13} \; & = \Lg_{7,1.5} \hspace{16cm} \\
[x_1,x_2] & = x_3; \; [x_1,x_3] = x_4; \; [x_1,x_4] = x_5; \; [x_1,x_5] = x_6; \\ 
[x_2,x_3] & =x_6; \; [x_2,x_5]=-x_7;\;  [x_3,x_4]=x_7.
\end{align*}
\begin{align*}
\Lg_{14} \; & = \Lg_{7,1.10} \hspace{16cm} \\
[x_1,x_2] & = x_3; \; [x_1,x_3] = x_4; \; [x_1,x_4] = x_6; \; 
[x_1,x_6] = x_7; \; 
[x_2,x_3] =x_5; \; [x_2,x_5]=x_7.
\end{align*}
\begin{align*}
\Lg_{15} \; & = \Lg_{7,1.11} \hspace{16cm} \\
[x_1,x_2] & = x_4; \; [x_1,x_4] = x_5; \; [x_1,x_5] = x_6; \; 
[x_1,x_6] = x_7; \; \\
[x_2,x_3] & =x_6; \; [x_2,x_4] =x_6; \; [x_2,x_5]=x_7;\; [x_3,x_4]=-x_7.
\end{align*}
\begin{align*}
\Lg_{16} \; & = \Lg_{7,1.14} \hspace{16cm} \\
[x_1,x_2] & = x_3; \; [x_1,x_3] = x_4; \; [x_1,x_4] = x_5+x_6; \;
[x_2,x_3] =x_5; \; [x_2,x_5]=-x_7;\; [x_3,x_4]=x_7.
\end{align*}
\begin{align*}
\Lg_{17} \; & = \Lg_{7,1.17} \hspace{16cm} \\
[x_1,x_2] & = x_3; \; [x_1,x_3] = x_4; \; [x_1,x_4] = x_6; \; 
[x_1,x_6] = x_7; \; \\
[x_2,x_3] & =x_5; \; [x_2,x_5] =x_6; \; [x_2,x_6]=x_7;\; [x_3,x_4]=-x_7; \;
[x_3,x_5]=x_7.
\end{align*}
\begin{align*}
\Lg_{18} \; & = \Lg_{7,1.21} \hspace{16cm} \\
[x_1,x_2] & = x_4; \; [x_1,x_4] = x_5; \; [x_1,x_5] = x_6; \\
[x_2,x_3] & =x_6; \; [x_2,x_4] =x_6; \; [x_2,x_6]=x_7;\; [x_4,x_5]=-x_7.
\end{align*}
\begin{align*}
\Lg_{D} \; & = \Lg_{7,1.4} \hspace{16cm} \\
[x_1,x_i] & = x_{i+1}; \; 2\le i\le 6 ;\quad
[x_2,x_3] =x_6; \; [x_2,x_4] =x_7.
\end{align*}
\begin{align*}
\Lg_{B} \; & = \Lg_{7,1.6} \hspace{16cm} \\
[x_1,x_i] & = x_{i+1}; \; 2\le i\le 6 ;\quad
[x_2,x_3] =x_7.
\end{align*}
\begin{align*}
\Lg_{19} \; & = \Lg_{7,1.01(ii)} \hspace{16cm} \\
[x_1,x_2] & = x_4; \; [x_1,x_4] = x_5; \; [x_1,x_5] = x_6;  \; [x_1,x_6] = x_7; \\
[x_2,x_3] & =x_6+x_7; \; [x_3,x_4]=-x_7.
\end{align*}
\begin{align*}
\Lg_{20} \; & = \Lg_{7,1.12} \hspace{16cm} \\
[x_1,x_2] & = x_4; \; [x_1,x_4] = x_5; \; [x_1,x_5] = x_6;  \; [x_1,x_6] = x_7; \\
[x_2,x_3] & =x_7; \; [x_2,x_4] =x_6; \; [x_2,x_5]=x_7.
\end{align*}
\begin{align*}
\Lg_{21} \; & = \Lg_{7,1.13} \hspace{16cm} \\
[x_1,x_2] & = x_3; \; [x_1,x_3] = x_4; \; [x_1,x_4] = x_6;  \; [x_1,x_5] = x_7; 
 \; [x_1,x_6] = x_7;\\
[x_2,x_3] & =x_5; \; [x_2,x_4] =x_7.
\end{align*}
\begin{align*}
\Lg_{22} \; & = \Lg_{7,2.4} \hspace{16cm} \\
[x_1,x_2] & = x_3; \; [x_1,x_3] = x_4; \; [x_1,x_4] = x_5;  \; [x_1,x_5] = x_6; \;
[x_2,x_5] =-x_7; \; [x_3,x_4] =x_7.
\end{align*}
\begin{align*}
\Lg_{23} \; & = \Lg_{7,2.5} \hspace{16cm} \\
[x_1,x_2] & = x_3; \; [x_1,x_3] = x_4; \; [x_1,x_5] = x_6;  \; [x_1,x_6] = x_7;\\
[x_2,x_3] & =x_5; \; [x_2,x_4] =x_6; \;[x_3,x_4] =x_7.
\end{align*}
\begin{align*}
\Lg_{24} \; & = \Lg_{7,2.6} \hspace{16cm} \\
[x_1,x_2] & = x_3; \; [x_1,x_3] = x_4; \; [x_1,x_4] = x_5;  \;
[x_2,x_3] =x_6; \; [x_2,x_5] =x_7; \;[x_3,x_4] =-x_7.
\end{align*}
\begin{align*}
\Lg_{25} \; & = \Lg_{7,2.10} \hspace{16cm} \\
[x_1,x_2] & = x_3; \; [x_1,x_3] = x_4; \; [x_1,x_4] = x_5;  \; [x_1,x_6] = x_7; \;
[x_2,x_5] = x_7; \;[x_3,x_4] =-x_7.
\end{align*}
\begin{align*}
\Lg_{26} \; & = \Lg_{7,2.13} \hspace{16cm} \\
[x_1,x_2] & = x_4; \; [x_1,x_4] = x_5; \; [x_1,x_5] = x_6; \;
[x_2,x_3] = x_6; \;[x_2,x_6] =x_7;\; [x_4,x_5] =-x_7.
\end{align*}
\begin{align*}
\Lg_{27} \; & = \Lg_{7,2.14} \hspace{16cm} \\
[x_1,x_3] & = x_4; \; [x_1,x_4] = x_5; \; [x_1,x_5] = x_6; \; [x_1,x_6] = x_7; \\
[x_2,x_3] & = x_5; \;[x_2,x_4] =x_6;\; [x_2,x_5] =x_7.
\end{align*}
\begin{align*}
\Lg_{A} \; & = \Lg_{7,2.3} \hspace{16cm} \\
[x_1,x_i] & = x_{i+1}; \; 2\le i\le 6 
\end{align*}
\begin{align*}
\Lg_{28} \; & = \Lg_{7,1.15} \hspace{16cm} \\
[x_1,x_2] & = x_4; \; [x_1,x_4] = x_5; \; [x_1,x_5] = x_6; \; [x_1,x_6] = x_7; \;
[x_2,x_3] = x_7; \;[x_2,x_4] =x_7.
\end{align*}
\begin{align*}
\Lg_{29} \; & = \Lg_{7,2.7} \hspace{16cm} \\
[x_1,x_2] & = x_3; \; [x_1,x_3] = x_4; \; [x_1,x_4] = x_6; \; [x_1,x_6] = x_7; \;
[x_2,x_3] = x_5.
\end{align*}
\begin{align*}
\Lg_{30} \; & = \Lg_{7,2.15} \hspace{16cm} \\
[x_1,x_2] & = x_4; \; [x_1,x_4] = x_5; \; [x_1,x_5] = x_6; \; [x_1,x_6] = x_7; \;
[x_2,x_3] = x_6; \;[x_3,x_4] =-x_7.
\end{align*}
\begin{align*}
\Lg_{31} \; & = \Lg_{7,2.16} \hspace{16cm} \\
[x_1,x_2] & = x_4; \; [x_1,x_4] = x_5; \; [x_1,x_5] = x_6; \; [x_1,x_6] = x_7; \;
[x_2,x_3] = x_7.\\
\end{align*}

\begin{defi}
Let $\Lg$ be in $\CN_7(\C)$. We say that $\Lg$ admits a basis
of type $I$ if there is a basis $(x_1,\ldots, x_7)$ of $\Lg$
such that $[x_i,x_j]=0$ for all $1\le i,j\le 7$ with $i+j>7$.
\end{defi}

Only $6$ algebras of the above list do not admit a basis of type $I$.
For the other ones we have chosen such a basis. That means that we have
replaced the basis used in \cite{MAG} for the following algebras:
$\Lg_4,\; \Lg_8,\; \Lg_{16},\;\Lg_{24},\; \Lg_{25}$ and $\Lg_{27}$.
By an explicit computation the following lemma is easy to verify.

\begin{lem}
The only algebras of the above list which do not admit a basis
of type $I$ are $\Lg_{1}(\la),\;\Lg_{2}\; ,\Lg_{5}, \; \Lg_{17}, \; \Lg_{18}$
and $\Lg_{26}.$ 
\end{lem}

In the case of the above six algebras the computations for the
degenerations become very complicated and 
we will exclude these algebras from the study of degenerations.\\[0.2cm]
Let $h_i=\dim H^i(\Lg,\Lg)$ respectively  $b_i=\dim H^i(\Lg,\C)$
be the dimensions of the adjoint cohomo\-lo\-gy and the trivial cohomology.
Let $\al_1,\al_2$ be the complex roots of the polynomial
$x^2-x+1$ and $A=\{0,-2,1-\al_1,1-\al_2\}$.
The next table gives a summary of some invariants of our algebras.
Let $n(\Lg)$ respectively $s(\Lg)$ denote the nilpotency and solvability class
of $\Lg$. 
\vspace*{0.3cm}
\begin{center}
\begin{tabular}{c|c|c|c|c|c|c}
 $\Lg$ & $(h_0,h_1,h_2,h_3,\ldots ,h_7)$ &  $(b_1,b_2,b_3,\ldots ,b_7)$ 
& $n(\Lg)$ & $s(\Lg)$  & $\dim O(\Lg)$ & $\Lg/Z(\Lg)$ \\
\hline
$\Lg_I(\al), \al \not \in A$ & $(1,4,9,14,15,11,6,2)$ & $(2,3,4,4,3,2,1)$ 
& $6$ & $3$  & $39$ & $\Lg_{6,D}$\\ 
$\Lg_I(-2)$ & $(1,4,9,15,16,12,7,2)$ & $(2,4,5,5,4,2,1)$ & $6$ & $3$ & $39$ & $\Lg_{6,D}$\\
$\Lg_I(1-\al_i)$ & $(1,4,9,14,16,12,6,2)$ & $(2,3,5,5,3,2,1)$ & $6$ & $3$ 
& $39$ & $\Lg_{6,D}$\\ 
$\Lg_F$ & $(1,4,9,15,16,11,6,2)$  & $(2,3,4,4,3,2,1)$ & $6$ & $3$ & $39$ & $\Lg_{6,B}$ \\
$\Lg_H$ & $(1,4,10,15,15,11,6,2)$ & $(2,3,4,4,3,2,1)$ & $6$ & $2$ & $39$ & $\Lg_{6,D}$ \\
$\Lg_3$ & $(1,4,9,15,17,13,7,2)$ & $(2,3,5,5,3,2,1)$ & $5$ & $3$ & $39$ & $1346_C$ \\
$\Lg_4$ & $(1,4,10,17,18,13,7,2)$ & $(2,3,4,4,3,2,1)$ & $5$ & $3$ & $39$ & $2346$ \\
$\Lg_6$ & $(2,5,9,14,15,11,6,2)$ & $(2,3,4,4,3,2,1)$ & $5$ & $3$ & $39$ & $1235_B$ \\
$\Lg_7$ & $(1,4,11,16,16,15,10,3)$ & $(3,4,4,4,4,3,1)$ & $5$ & $2$ & $39$ & $1246$ \\
\hline
$\Lg_C$ & $(1,5,10,15,16,11,6,2)$ & $(2,3,4,4,3,2,1)$ & $6$ & $3$ & $38$ & $\Lg_{6,A}$ \\
$\Lg_G$ & $(1,5,11,15,15,11,6,2)$ & $(2,3,4,4,3,2,1)$ & $6$ & $2$ & $38$ & $\Lg_{6,D}$ \\
$\Lg_E$ & $(1,5,12,19,20,14,7,2)$ & $(2,4,6,6,4,2,1)$ & $6$ & $2$ & $38$ & $\Lg_{6,B}$ \\
$\Lg_8$ & $(1,5,13,17,16,15,10,3)$ & $(3,4,4,4,4,3,1)$ & $5$ & $2$ & $38$ & $1246$ \\
$\Lg_9$ & $(1,5,11,16,17,13,7,2)$ & $(2,3,5,5,3,2,1)$ & $5$ & $3$ & $38$ & $1346_C$ \\
$\Lg_{10}$ & $(1,5,11,17,18,13,7,2)$ & $(2,3,4,4,3,2,1)$ & $5$ & $3$ & $38$ & $2346$ \\
$\Lg_{11}$ & $(1,5,10,15,18,14,7,2)$ & $(2,3,6,6,3,2,1)$ & $5$ & $3$ & $38$ & $1346_C$ \\
$\Lg_{12}$ & $(1,5,12,18,19,16,10,3)$ & $(3,4,4,4,4,3,1)$ & $5$ & $3$ & $38$ & $1+1235_B$ \\
$\Lg_{13}$ & $(2,6,11,17,17,11,6,2)$ & $(2,3,4,4,3,2,1)$ & $5$ & $3$ & $38$ & $1235_A$ \\
$\Lg_{14}$ & $(1,5,12,19,21,16,8,2)$ & $(2,4,7,7,4,2,1)$ & $5$ & $2$ & $38$ & $2346$ \\
$\Lg_{15}$ & $(1,5,14,22,23,19,11,3)$ & $(3,5,6,6,5,3,1)$ & $5$ & $2$ & $38$ & $1346_B$ \\
$\Lg_{16}$ & $(2,6,11,17,18,13,7,2)$ & $(2,3,4,4,3,2,1)$ & $5$ & $3$ & $38$ & $1235_B$ \\
\hline
$\Lg_D$ & $(1,6,13,19,20,14,7,2)$ & $(2,4,6,6,4,2,1)$ & $6$ & $2$ & $37$ & $\Lg_{6,B}$ \\
$\Lg_B$ & $(1,6,15,23,22,14,7,2)$ & $(2,4,6,6,4,2,1)$ & $6$ & $2$ & $37$ & $\Lg_{6,A}$ \\
$\Lg_{19}$ & $(1,6,16,24,25,20,11,3)$ & $(3,5,7,7,5,3,1)$ & $5$ & $2$ & $37$ & $1346_A$ \\
$\Lg_{20}$ & $(1,6,16,26,28,21,11,3)$ & $(3,5,7,7,5,3,1)$ & $5$ & $2$ & $37$ & $1+1235_B$ \\
$\Lg_{21}$ & $(2,7,14,21,22,16,8,2)$ & $(2,4,7,7,4,2,1)$ & $5$ & $2$ & $37$ & $1235_A$ \\
$\Lg_{22}$ & $(2,7,12,17,17,12,7,2)$ & $(2,3,4,4,3,2,1)$ & $5$ & $3$ & $37$ & $1235_A$ \\
$\Lg_{23}$ & $(1,6,12,16,18,14,7,2)$ & $(2,3,6,6,3,2,1)$ & $5$ & $3$ & $37$ & $1346_C$ \\
$\Lg_{24}$ & $(2,7,12,17,18,13,7,2)$ & $(2,3,4,4,3,2,1)$ & $5$ & $3$ & $37$ & $1235_A$ \\
$\Lg_{25}$ & $(1,6,14,19,19,16,10,3)$ & $(3,4,4,4,4,3,1)$ & $5$ & $3$ & $37$ & $1+1235_A$ \\
$\Lg_{27}$ & $(1,6,14,17,16,15,10,3)$ & $(3,4,4,4,4,3,1)$ & $5$ & $2$ & $37$ & $1246$ \\
\hline
$\Lg_A$ & $(1,7,17,25,23,14,7,2)$ & $(2,4,6,6,4,2,1)$ & $6$ & $2$ & $36$ & $\Lg_{6,A}$ \\
$\Lg_{28}$ & $(1,7,18,27,28,21,11,3)$ & $(3,5,7,7,5,3,1)$ & $5$ & $2$ & $36$ & $1+1235_A$ \\
$\Lg_{29}$ & $(2,8,16,25,25,16,8,2)$ & $(2,4,7,7,4,2,1)$ & $5$ & $2$ & $36$ & $1235_A$ \\
$\Lg_{30}$ & $(1,7,18,26,26,20,11,3)$ & $(3,5,7,7,5,3,1)$ & $5$ & $2$ & $36$ & $1346_A$ \\
\hline
$\Lg_{31}$ & $(1,8,20,28,28,21,11,3)$ & $(3,5,7,7,5,3,1)$ & $5$ & $2$ & $35$ & $1+1235_A$\\
\end{tabular}
\end{center}   
\vspace*{0.5cm}
We have $b_0=1$ for all these algebras. Therefore we have omitted it in
the list. The central quotients $\Lg/Z(\Lg)$ are nilpotent Lie algebras of 
dimension $5$ and $6$. We have used the notation from \cite{SEE}. 
Note that $\Lg_I(0)=\Lg_G$.\\[0.2cm]
We divide the classification of all degenerations according to the
orbit dimensions. If the orbit dimension of $O(\la)$ is smaller or equal
than the dimension of $O(\mu)$, then $\la$ cannot degenerate to 
$\mu$.

\begin{prop}
All non-trivial degenerations between algebras of the above table with
orbit dimension $38$ and $39$ are given as follows:
\vspace*{0.5cm}
\begin{center}
\begin{tabular}{c|c|c|c|c|c|c|c|c|c|c|c|c}
$ \ra_{\rm deg}$ & $\Lg_C$ & $\Lg_G$ & $\Lg_E$ & $\Lg_8$ & $\Lg_9$ & $\Lg_{10}$ & $\Lg_{11}$ &
$\Lg_{12}$ & $\Lg_{13}$ & $\Lg_{14}$ & $\Lg_{15}$ & $\Lg_{16}$ \\
\hline
$\Lg_I(\al),\al \neq 1$ & B & B & $\ch$ & B & B & B & B & B & B & $\ch$ & $\ch$ & B \\
$\Lg_I(1)$ & B & I & $\ch$ & I & $\ch$ & B & I & I & B & I & B & I  \\
$\Lg_F$ & $\ch$ & z &  $\ch$ & $z_3$ & z &  $\ch$ & z & I &  $\ch$ & I & z & I \\
$\Lg_H$ & $z_2$ & $\ch$ & $\ch$ & $\ch$ & s & s & s & s & s & $\ch$ & B & s \\
$\Lg_3$ & n & n & n & $z_3$ & $\ch$ & $b_3$ & $\ch$ & $b_3$ & $h_5$ & $\ch$ & $\ch$ & $b_3$ \\
$\Lg_4$ & n & n & n & $z_3$ & $h_3$ & $\ch$ & $z_2$ & $\ch$ & $h_5$ & $\ch$ & $z$ & $\ch$ \\
$\Lg_6$ & n & n & n & $h_0$ & $h_0$ & $h_0$ & $h_0$ & $h_0$ & $\ch$ & $h_0$ & $h_0$ &$\ch$ \\
$\Lg_7$ & n & n & n & $\ch$ & $h_5$ & $h_5$ & $h_2$ & s & $h_5$ & $h_6$ & $\ch$ & $h_5$ \\
\end{tabular}
\end{center}   
\vspace*{0.5cm}
\end{prop}

\begin{proof}
The checkmark denotes that there is a degeneration
$\la \ra_{\rm deg} \mu$. The other symbols stand for the reason
why such a degeneration is impossible. In general there is more
than just one reason for a non-degeneration. However we have
written down only one in the table. \\
The symbol $z$ denotes the fact
that $\la$ cannot degenerate to $\mu$ if the central quotients do not degenerate
to each other.  Here we use the result,
that if a nilpotent Lie algebra $\Lg$ degenerates to $\Lh$, then the
central quotient $\Lg/Z(\Lg)$ degenerates to $\Lh/Z(\Lh) \oplus \C^d$,
where $\C^d$ is an abelian Lie algebra of dimension $d=\dim Z(\Lh)-\dim Z(\Lg)$,
see \cite{SEE1}. For example, $\Lg_F$ cannot degenerate to $\Lg_G$ since
$12346_B$ does not degenerate to $12346_D$ in dimension $6$. \\
The symbols $h_i$ denote the fact, that $\la$ cannot degenerate to $\mu$
if $h_i(\la)=\dim H^i(\la,\la) > \dim H^i(\mu,\mu)=h_i(\mu)$ for some $i$. 
As an example consider $\Lg_7$ and $\Lg_9$ where $h_5(\Lg_7)=15$ and
$h_5(\Lg_9)=13$. Similarly $b_i= \dim H^i(\la)$ and $z_i=\dim Z^i(\la,\la)$ are used.
Note that $z_3(\Lg_8)=113,\; z_3(\Lg_F)=z_3(\Lg_3)=114,\; z_3(\Lg_4)=115$
and $z_2(\Lg_H)=49,\; z_2(\Lg_C)=48$. \\
The symbols $n$ and $s$ stand for nilpotency and solvability class of $\Lg$.
If $\la \ra_{\rm deg} \mu$, then $n(\la)\ge n(\mu)$ and $s(\la)\ge s(\mu)$.\\
The symbol $I$ denotes the following fact. If $\la$ degenerates to $\mu$
and $\la$ is represented by a structure, which lies in a $B$-stable subset
$\CR$ of $\CN_7(\C)$ for some Borel subgroup $B$ in $G=GL(7,\C)$, then $\mu$ must
also be represented by a structure in $\CR$. Let $\CR$ be defined by the
property that $\Lg$ possesses an ideal $I$ of codimension $1$ such that
\begin{align*}
[\Lg,\Lg] & \subseteq I \\
[\Lg,[\Lg,\Lg]] & = 0
\end{align*}   
It is obvious that $\Lg_I(1)$ and $\Lg$ admit such an ideal: $I=\s\{x_2,\ldots,x_7\}$.
On the other hand, $\Lg_G,\;\Lg_8,\; \Lg_{11},\; \Lg_{12},\; \Lg_{14}$ and
$\Lg_{16}$ do not admit such an ideal. Hence there is no degeneration from
$\Lg_I(1)$ and $\Lg_F$ to these algebras. Let us show, as an example, why $\Lg=\Lg_G$
does not admit such an ideal. Because of $[\Lg,\Lg] \subseteq I$ we would have
$I=\s\{y,x_3,x_4,x_5,x_6,x_7\}$ with $y=\al x_1+\be x_2$. Then
\begin{align*}
[y,[y,x_3]] & = [\al x_1+\be x_2,\al x_4+\be x_5]=\al^2 x_5+2\al\be x_6+\be^2 x_7
\end{align*}
Hence $[\Lg,[\Lg,\Lg]] = 0$ would imply $\al=\be=0$ and $y=0$. This contradicts
$\dim I=6$.\\
The symbol $B$ stands for the following argument. Let $B$ be the Borel subgroup
of $G=GL(7,\C)$ consisting of invertible lower-triangular matrices.
Then we have $\ov{G\cdot \mu}=G\cdot \ov{B\cdot \mu}$ for all $\mu \in \CN_7(\C)$,
see \cite{GRH}.
If we can show that $\mu$ is not isomorphic to any algebra contained in the closure of 
the $B$-orbit of $\la$, then $\la$ cannot degenerate to $\mu$.
Consider the $B$-orbit of $\Lg_I(\al)$. It consists of algebras $\Lg(\al_1,\ldots, \al_{22})$
with Lie brackets
\begin{align*}
[x_i,x_j] & = \sum_{k=i+j}^7 \al_{i,j}^k x_k
\end{align*} 
where $\al_1=\al_{1,2}^3,\; \al_2=\al_{1,2}^4,\ldots , \al_{22}=\al_{3,4}^7$.
The algebra $\Lg_C$ is isomorphic to  $\Lg(\al_1,\ldots, \al_{22})$ if and only
if certain conditions on the $\al_i$ are satisfied. ( Necessary conditions in this case
are $a_{16}=a_{17}=a_{19}=0$ and $\al_1,\al_6,\al_{10},a_{13},a_{15},a_{21},a_{22}$ non-zero ). 
However, it is easy to see by an explicit computation that in the closure of the 
$B$-orbit of $\Lg_I(\al)$ there is no such algebra satisfying these conditions.
Hence $\Lg_I(\al)$ does not degenerate to $\Lg_C$. (For $\al \in A$ this follows also
from the adjoint cohomology $h_5$). Later we will see that $\Lg_C  \ra_{\rm deg} \Lg_{25}$ and
$\Lg_I(\al)$ cannot degenerate to $\Lg_{25}$. By transitivity it follows again that
$\Lg_I(\al)$ does not degenerate to $\Lg_C$. \\
In case there is a checkmark in the table we have found a degeneration
$\la \ra_{\rm deg}\mu$ by explicitely constructing a $g_t\in GL(7,\C(t))$ such that   
$\mu\cong\lim_{t \to 0} g_t\cdot\la$.
The degenerations of $\Lg_F$ are as follows. $\Lg_F \ra_{\rm deg}\Lg_E$
can be realized by
$$ g_t^{-1}=\begin{pmatrix} 
t & 0 & 0 & 0 & 0 & 0 & 0\\ 
0 & t^3 & 0 & 0 & 0 & 0 & 0\\ 
0 & 0 & t^4 & 0 & 0 & 0 & 0\\
0 & t^4/2 & 0 & t^5 & 0 & 0 & 0\\
0 & 0 & t^5/2 & 0 & t^6 & 0 & 0 \\
0 & 0 & 0 & t^6/2 & 0 & t^7 & 0\\
0 & 0 & 0 & 0 & t^7/2 & 0 & t^8
\end{pmatrix}$$
The other ones are realized by diagonal matrices.
\begin{align*}
\Lg_F \ra_{\rm deg}\Lg_C,\quad g_t^{-1}& =\diag (t^{-1},t^{-2},t^{-3},t^{-4},
t^{-5},t^{-6},t^{-7})\\  
\Lg_F \ra_{\rm deg}\Lg_D,\quad g_t^{-1}& =\diag (t,t^3,t^4,t^5,t^6,t^7,t^8)\\
\Lg_F \ra_{\rm deg}\Lg_{13},\quad g_t^{-1}& =\diag (t^{-1},t^{-3},t^{-4}
,t^{-5},t^{-6},t^{-7},-t^{-9})  
\end{align*}
Here $\Lg_F$ degenerates to the algebras with the Lie brackets exactly as given
in the list. In general however, if $\la \ra_{\rm deg} \mu$, then $\mu$
is only isomorphic to the algebra given in our list.  
We have a complete list of all degeneration matrices. It is however too long
to be given here.
\end{proof}

\begin{prop}
All non-trivial degenerations from algebras of 
orbit dimension $39$ to algebras of orbit dimension $35,36,37$ are given as follows:
\vspace*{0.5cm}
\begin{center}
\begin{tabular}{c|c|c|c|c|c|c|c|c|c|c|c|c|c|c|c}
$ \ra_{\rm deg}$ & $\Lg_D$ & $\Lg_B$ & $\Lg_{19}$ & $\Lg_{20}$ & $\Lg_{21}$ & $\Lg_{22}$ 
& $\Lg_{23}$ & $\Lg_{24}$ & $\Lg_{25}$ & $\Lg_{27}$  & $\Lg_{A}$  & $\Lg_{28}$  & $\Lg_{29}$
 & $\Lg_{30}$  & $\Lg_{31}$ \\
\hline
$\Lg_I(\al),\al \neq 1$ & $\ch$ & $\ch$ & $\ch$ & $\ch$ &  $\ch$ & B & B & B & B & B 
& $\ch$ & $\ch$ & $\ch$ & $\ch$ & $\ch$ \\
$\Lg_I(1)$ & $\ch$ & $\ch$ & $\ch$ & $\ch$ & $\ch$ & B & $\ch$ & B & B & B  
& $\ch$ & $\ch$ & $\ch$ & $\ch$ & $\ch$ \\
$\Lg_F$ & $\ch$ & $\ch$ & $\ch$ & $\ch$ & $\ch$ &  $\ch$ & z & $\ch$ &$\ch$ & $z_3$
& $\ch$ & $\ch$ & $\ch$ & $\ch$ & $\ch$  \\
$\Lg_H$ & $\ch$ & $\ch$ & $\ch$ & $\ch$ & $\ch$ & s & s & s & s & $\ch$
& $\ch$ & $\ch$ & $\ch$ & $\ch$ & $\ch$  \\
$\Lg_3$ & n & n & $\ch$ & $\ch$ & $\ch$ & $h_5$ & $\ch$ & $b_3$ & $b_3$ & $h_4$ 
& n & $\ch$ & $\ch$ & $\ch$ & $\ch$ \\
$\Lg_4$ & n & n & z & $\ch$ & $\ch$ & $h_4$ & $h_3$ & $\ch$ & $\ch$ & $h_4$
& n & $\ch$ & $\ch$ & z & $\ch$ \\
$\Lg_6$ & n & n & $h_0$ & $h_0$ & $\ch$ & $\ch$ & $h_0$ & $\ch$ & $h_0$ & $h_0$
& n & $h_0$ & $\ch$ & $h_0$ & $h_0$ \\
$\Lg_7$ & n & n & $\ch$ & $\ch$ & $h_6$ & $h_5$ & $h_5$ & $h_5$ & s & $\ch$
 & n & $\ch$ & $h_6$ & $\ch$ & $\ch$ 
\end{tabular}
\end{center}   
\vspace*{0.5cm}
\end{prop}

\begin{prop}
All non-trivial degenerations from algebras of 
orbit dimension $38$ to algebras of orbit dimension $35,36,37$ are given as follows:
\vspace*{0.5cm}
\begin{center}
\begin{tabular}{c|c|c|c|c|c|c|c|c|c|c|c|c|c|c|c}
$ \ra_{\rm deg}$ & $\Lg_D$ & $\Lg_B$ & $\Lg_{19}$ & $\Lg_{20}$ & $\Lg_{21}$ & $\Lg_{22}$ 
& $\Lg_{23}$ & $\Lg_{24}$ & $\Lg_{25}$ & $\Lg_{27}$  & $\Lg_{A}$  & $\Lg_{28}$  & $\Lg_{29}$
 & $\Lg_{30}$  & $\Lg_{31}$ \\
\hline
$\Lg_C$ & z & $\ch$ & z & z & B & $\ch$ & z & B & $\ch$ & $z_3$ 
& $\ch$ & $\ch$ & $\ch$ & z & $\ch$ \\
$\Lg_G$ & B & $\ch$ & B & $\ch$ & $\ch$ & s & s & s & s & B  
& $\ch$ & $\ch$ & $\ch$ & B & $\ch$ \\
$\Lg_E$ & $\ch$ & $\ch$ & $\ch$ & $\ch$ & $\ch$ & $h_3$ & $h_3$ & $h_3$ &$h_4$ & $h_3$
& $\ch$ & $\ch$ & $\ch$ & $\ch$ & $\ch$  \\
$\Lg_{8}$ & n & n & $\ch$ & $\ch$ & $h_6$ & $h_2$ & $h_2$ & $h_2$ & s & $\ch$
& n & $\ch$ & $h_6$ & $\ch$ & $\ch$  \\
$\Lg_{9}$ & n & n & $\ch$ & $\ch$ & $\ch$ & $h_5$ & $\ch$ & $b_3$ & $b_3$ & $h_4$ 
& n & $\ch$ & $\ch$ & $\ch$ & $\ch$ \\
$\Lg_{10}$ & n & n & z & $\ch$ & $\ch$ & $h_4$ & $h_3$ & $\ch$ & $\ch$ & $h_4$
& n & $\ch$ & $\ch$ & B & $\ch$ \\
$\Lg_{11}$ & n & n & $\ch$ & $\ch$ & $\ch$ & $h_4$ & $\ch$ & $h_5$ & $b_3$ & $h_4$
& n & $\ch$ & $\ch$ & $\ch$ & $\ch$ \\
$\Lg_{12}$ & n & n & z & $\ch$ & $h_6$ & $h_3$ & $h_3$ & $h_3$ & $\ch$ & $h_3$
& n & $\ch$ & $h_6$ & z & $\ch$ \\
$\Lg_{13}$ & n & n & $h_0$ & $h_0$ & B & $\ch$ & $h_3$ & $\ch$ & $h_0$ & $h_4$
& n & $h_0$ & $\ch$ & $h_0$ & $h_0$ \\
$\Lg_{14}$ & n & n & z & $\ch$ & $\ch$ & $h_3$ & $h_3$ & $h_3$ & $h_4$ & $h_3$
& n & $\ch$ & $\ch$ & z & $\ch$ \\
$\Lg_{15}$ & n & n & $\ch$ & $\ch$ & $h_3$ & $h_2$ & $h_2$ & $h_2$ & $h_3$ & $h_3$
& n & $\ch$ & $h_5$ & $\ch$ & $\ch$ \\
$\Lg_{16}$ & n & n & $h_0$ & $h_0$ & $\ch$ & $h_4$ & $h_3$ & $\ch$ & $h_0$ & $h_4$
& n & $h_0$ & $\ch$ & $h_0$ & $h_0$ 
\end{tabular}
\end{center}   
\vspace*{0.5cm}
\end{prop}

\begin{prop}
All non-trivial degenerations between algebras of 
orbit dimension $37$ and algebras of orbit dimension $35,36$ are given as follows:
\vspace*{0.5cm}
\begin{center}
\begin{tabular}{c|c|c|c|c|c}
$ \ra_{\rm deg}$ & $\Lg_A$ & $\Lg_{28}$ & $\Lg_{29}$ & $\Lg_{30}$ & $\Lg_{31}$ \\ 
\hline
$\Lg_D$ & $\ch$ & B & B & $\ch$ & B \\ 
$\Lg_B$ & $\ch$ & $\ch$ & $\ch$ & z & $\ch$ \\  
$\Lg_{19}$ & n & $\ch$ & $h_5$ & $\ch$ & $\ch$  \\
$\Lg_{20}$ & n & $\ch$ & $h_5$ & z & $\ch$ \\
$\Lg_{21}$ & n & $h_0$ & $\ch$ & $h_0$ & $h_0$ \\
$\Lg_{22}$ & n & $h_0$ & $\ch$ & $h_0$ & $h_0$ \\
$\Lg_{23}$ & n & B & $\ch$ & $\ch$ & B \\
$\Lg_{24}$ & n & $h_0$ & $\ch$ & $h_0$ & $h_0$ \\
$\Lg_{25}$ & n & $\ch$ & $h_6$ & z & $\ch$ \\
$\Lg_{27}$ & n & $\ch$ & $h_6$ & B & $\ch$  
\end{tabular}
\quad \quad
\begin{tabular}{c|c}
$ \ra_{\rm deg}$ & $\Lg_{31}$ \\ 
\hline
$\Lg_{28}$ & $\ch$ \\
$\Lg_{29}$ & $h_0$ \\
$\Lg_{30}$ & B \\
$\Lg_{31}$ & B 
\end{tabular}
\end{center}   
\vspace*{0.5cm}
\end{prop}

\begin{proof}
The use of transitivity for degenerations is very helpful. As an example,
we obtain all possible degenerations of $\Lg_F$ to algebras of orbit dimension
$35,36,37$ by the degenerations $\Lg_F\ra_{\rm deg} \Lg_C$, $\Lg_F\ra_{\rm deg} \Lg_E$,
$\Lg_F\ra_{\rm deg} \Lg_{10}$
and the degenerations of the algebras $\Lg_C$, $\Lg_E$ and  $\Lg_{10}$:
If we degenerate $\Lg_E$ via $ g_t^{-1}=\diag (t^2,t^7,t^9,t^{11},t^{13},t^{15},t^{16})$
then we obtain a Lie algebra with defining brackets $[y_1,y_i]=y_{i+1},\; 2\le i\le 5$
and $[y_2,y_3]=y_7$. It is isomorphic to $\Lg_{29}$ by setting $y_5=x_6$, 
$y_6=x_5$ and $y_i=x_i$ otherwise. Similarly we obtain:
\begin{align*}
\Lg_E \ra_{\rm deg}\Lg_D,\quad g_t^{-1}& =\diag (t^{-1},t^{-3},t^{-4},t^{-5},
t^{-6},t^{-7},t^{-8})\\  
\Lg_E \ra_{\rm deg}\Lg_B,\quad g_t^{-1}& =\diag (t,t^4,t^5,t^6,t^7,t^8,t^9)\\
\Lg_E \ra_{\rm deg}\Lg_A,\quad g_t^{-1}& =\diag (t,t^5,t^6,t^7,t^8,t^9,t^{10})\\
\Lg_E \ra_{\rm deg}\Lg_{19},\quad g_t^{-1}& =\diag (1,1,t^{-1},t^{-1},t^{-1},t^{-1},
t^{-1})\\  
\Lg_E \ra_{\rm deg}\Lg_{30},\quad g_t^{-1}& =\diag (t^{-1},t^{-3},t^{-5},t^{-6},
t^{-7},t^{-8},t^{-9})\\   
\Lg_E \ra_{\rm deg}\Lg_{31},\quad g_t^{-1}& =\diag (t,t^{4},t^{5},1,
t^{2},t^{3},t^{4})\\
\Lg_C \ra_{\rm deg}\Lg_{22},\quad g_t^{-1}& =\diag (1,t^{-1},t^{-1},t^{-1},
t^{-1},t^{-1},-t^{-2})\\  
\Lg_C \ra_{\rm deg}\Lg_{25},\quad g_t^{-1}& =\diag (1,t^{-1},t^{-1},t^{-1},
t^{-1},t^{-2},t^{-2})  
\end{align*}
Furthermore $\Lg_E$ degenerates to $\Lg_{20}$ respectively to $\Lg_{21}$ by

$$ g_t^{-1}=\begin{pmatrix} 
-t & 0    & 0   & 0   & 0   & 0   & 0\\ 
0  & -t^2 & 0   & 0   & 0   & 0   & 0\\ 
0  & 0    & t^3 & 0   & 0   & 0   & 0\\
0  & t^2  & 0   &-t^4 & 0   & 0   & 0\\
0  & 0    & -t^3& 0   & t^4 & 0   & 0 \\
0  & 0    & 0   & t^4 & t^4 &-t^5 & 0\\
0  & 0    & 0   & 0   & 0   &-t^5 & t^6
\end{pmatrix},\quad  
g_t^{-1}=\begin{pmatrix} 
-t & 0    & 0   & 0   & 0   & 0   & 0\\ 
0  & -t^2 & 0   & 0   & 0   & 0   & 0\\ 
0  & 0    & t^3 & 0   & 0   & 0   & 0\\
0  & 0    & 0   &-t^4 & 0   & 0   & 0\\
0  & 0    & 0   & 0   & t^5 & 0   & 0 \\
0  & 0    & 0   & 0   & 0   &-t^5 & 0\\
0  & 0    & 0   & 0   & 0   &-t^5 & t^6
\end{pmatrix}
$$

The degenerations $\Lg_E \ra_{\rm deg}\Lg_{28}$ and $\Lg_{10} \ra_{\rm deg}\Lg_{24}$
can be realized by

$$ g_t^{-1}=\begin{pmatrix} 
t  & 0    & 0   & 0   & 0   & 0   & 0\\ 
0  & t^3  & 0   & 0   & 0   & 0   & 0\\ 
0  & 0    & t^3 & 0   & 0   & 0   & 0\\
0  &-t    &-t^2 & t^3 & 0   & 0   & 0\\
0  & 0    & 0   &-t^2 & t^4 & 0   & 0 \\
0  & 0    & 0   & 0   &-t^3 & t^5 & 0\\
0  & 0    & 0   & 0   & 0   &-t^4 & t^6
\end{pmatrix},\quad  
g_t^{-1}=\begin{pmatrix} 
t^{-1} & 0  & 0     & 0     & 0     & 0   & 0\\ 
-t^2   & 1  & 0     & 0     & 0     & 0   & 0\\ 
0      & 0  & t^{-1}& 0     & 0     & 0   & 0\\
0      & 0  & 0     & t^{-2}& 0     & 0   & 0\\
0      & 0  & 0     & 0     & t^{-3}& 0   & 0 \\
0      & 0  & 0     &-t     & 1   & t^{-1}& 0\\
0      & 0  & 0     & 0     &-2   & 0     & t^{-3}
\end{pmatrix}
$$

This shows that $\Lg_F$ degenerates to all algebras of orbit dimension
$35,36,37$ except for $\Lg_{23},\Lg_{27}$.
Transitivity is also useful for showing non-degenerations. Since $\Lg_I(\al)$
does not degenerate to $\Lg_{24}$ and $\Lg_{10},\,\Lg_{13},\,\Lg_{16}\ra_{\rm deg}
\Lg_{24}$ we conclude that  $\Lg_I(\al)$ cannot degenerate to $\Lg_{10},
\,\Lg_{13},\,\Lg_{16}$. Since $\Lg_I(\al)$ does not degenerate to $\Lg_{25}$, and
$\Lg_C,\Lg_{12}\ra_{\rm deg}\Lg_{25}$ it follows that $\Lg_I(\al)$ cannot degenerate
to $\Lg_C,\Lg_{12}$. Similarly we see that $\Lg_I(\al)$ cannot degenerate to
$\Lg_8$, and for $\al \ne 0$ not to $\Lg_9,\Lg_{11}$.
\end{proof}

\section{The varieties $\CL_n(k)$}

It is already quite interesting to investigate the varieties
$\CL_n(k)$ and the orbit closures over the complex numbers
in small dimensions. 
For $n=2$ we have 
$$\CL_2(\C)=\ov{O(\Lr_2(\C))}=O(\Lr_2(\C)) \cup O(\C^2)$$ 
where $\Lr_2(\C)$ is the non-abelian algebra. The only non-trivial
degeneration is given by $\Lr_2(\C)\ra_{\rm deg} \C^2$.
The orbit of $\Lr_2(\C)$ is open. There is no Lie algebra law
degenerating to $\Lr_2(\C)$ in $\CL_2(\C)$.\\
The variety $\CL_3(\C)$ is the union of two irreducible components 
$\CC_1$ and $\CC_2$. 
The component $\CC_1$ consists of the Lie algebras of trace zero,
i.e., where the linear form  ${\rm \mbox{tr ad}}(x)$ vanishes:
$$\CC_1=\ov{O(\Ls\Ll_2(\C))}=O(\Ls\Ll_2(\C)) \cup O(\Lr_{3,-1}(\C))
\cup O(\Ln_3(\C)) \cup O(\C^3)$$
The classification of all orbits and their orbit closures in
$\CL_3(\C)$ is given as follows:
\vspace*{0.5cm}
\begin{center}
\begin{tabular}{c|c|c}
 $\Lg$ & Lie brackets & $\ov{O(\Lg)}$  \\
\hline
$\C^3$ & $-$ & $\C^3$ \\
$\Ln_3(\C)$ & $[e_1,e_2]=e_3$ & $\Ln_3(\C),\;\C^3$ \\
$\Lr_2(\C) \oplus  \C$ & $[e_1,e_2]=e_2$ & $\Lr_2(\C)\oplus\C,\;\Ln_3(\C),\;\C^3$ \\
$\Lr_3(\C)$ & $[e_1,e_2]=e_2,\, [e_1,e_3]=e_2+ e_3$ & $\Ln_3(\C),\;\C^3$ \\ 
$\Lr_{3,\al}(\C)$ & $[e_1,e_2]=e_2, \,[e_1,e_3]=\al e_3,\, \abs{\al}< 1$
& $\Lr_{3,\al}(\C),\;\Ln_3(\C),\;\C^3 $\\
$\Lr_{3,-1}(\C)$ & $[e_1,e_2]=e_2, \,[e_1,e_3]=-e_3$
& $\Lr_{3,-1}(\C),\;\Ln_3(\C),\;\C^3 $\\
$\Lr_{3,1}(\C)$ & $[e_1,e_2]=e_2, \,[e_1,e_3]=e_3$ & $\Lr_{3,1}(\C),\;\C^3$ \\
$\Ls \Ll_2 (\C)$ & $[e_1,e_2]=e_3, [e_1,e_3]=-2 e_1,[e_2,e_3]=2 e_2$ & 
$\Ls \Ll_2 (\C),\;\Lr_{3,-1}(\C),\;\Ln_3(\C),\;\C^3$
\end{tabular}
\end{center}   
\vspace*{0.5cm}
The component $\CC_2$ consists of the solvable Lie algebras:
$$\CC_2=\CR_3(\C)=\cup_{\al}O(\Lr_{3,\al}(\C)) \cup O(\Lr_3(\C))\cup
O(\Lr_2(\C)\oplus \C)\cup O(\Ln_3(\C) \cup O(\C^3)$$
We have $\CC_1\cap \CC_2=\ov{O(\Lr_{3,-1}(\C))}$ and
$\dim \CC_1=\dim \CC_2=6$. 
The following diagram shows all essential degenerations (that is,
all the other degenerations are combinations of these) in 
$\CL_3(\C)$:

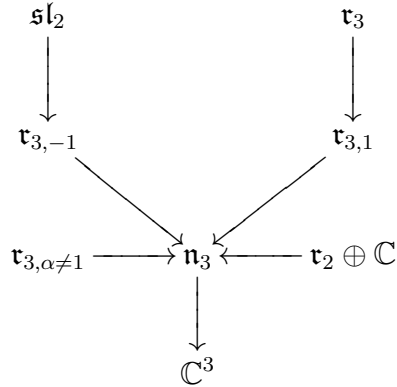
$$
\Dg
\Ls \Ll_2 & \phantom{--} & \phantom{--} & \phantom{--} & \Lr_3 \\
\dTo      &              &              &              & \dTo \\
\Lr_{3,-1}\aTo (2,-2)  &      &       &      &\Lr_{3,1}\aTo (-2,-2)  & \\
 & & & & & \\
\Lr_{3,\alpha \ne 1}   & \rTo & \Ln_3 & \lTo & \Lr_2\oplus \C \\
 & & \dTo & & \\
 & & \C^3 & & \\
\endDg
$$
In dimension $4$ the results become much more complicated.
\begin{prop}
The variety $\CL_4(\C)$ is the union of $4$ irreducible components 
$\CC_i$, $i=1,\ldots,4$ as follows:
\begin{align*}
\CC_1 & = \ov{O(\Ls\Ll_2(\C)\oplus \C)} \\
\CC_2 & = \ov{O(\Lr_2(\C)\oplus \Lr_2(\C))}\\
\CC_3 & = \ov{\cup_{\al,\be}O(\Lg_4(\al,\be))}\\
\CC_4 & = \ov{\cup_{\al}O(\Lg_5(\al))}
\end{align*}
\end{prop}
The components are of dimension $12$, i.e., $\dim \CC_i=12$.
The number of open orbits equals $2$; indeed, the Lie algebras
$\Ls\Ll_2(\C)\oplus \C$ and $\Lr_2(\C)\oplus \Lr_2(\C)$
are rigid.\\
The classification of all orbits in dimension $4$ is given in the
following table:
\vspace*{0.5cm}
\begin{center}
\begin{tabular}{c|c}
$\Lg$ & Lie brackets \\
\hline
$\C^4$ &  \\
$\Ln_3(\C)\oplus \C$ & $[e_1,e_2]=e_3$ \\
$\Ln_4(\C)$ & $[e_1,e_2]=e_3,\, [e_1,e_3]=e_4$ \\     
$\Lr_2(\C) \oplus  \C^2$ & $[e_1,e_2]=e_2$ \\
$\Lr_2(\C) \oplus \Lr_2(\C)$ & $[e_1,e_2]=e_2, \,[e_3,e_4]=e_4$ \\  
$\Ls\Ll_2(\C) \oplus \C$ & $[e_1,e_2]=e_2, \,[e_1,e_3]=-e_3,\,[e_2,e_3]=e_1$\\
$\Lg_1$ & $[e_1,e_2]=e_2, \,[e_1,e_3]=e_3,\,[e_1,e_4]=e_4$\\
$\Lg_2(\al)$ & $[e_1,e_2]=e_2, \,[e_1,e_3]=e_3,\,[e_1,e_4]=e_3+\al e_4$\\
$\Lg_3$ & $[e_1,e_2]=e_2, \,[e_1,e_3]=e_3,\,[e_1,e_4]=2e_4,\,[e_2,e_3]=e_4$\\
$\Lg_4(\al,\be)$ & $[e_1,e_2]=e_2, \,[e_1,e_3]=e_2+\al e_3,\,[e_1,e_4]=e_3+\be
e_4$\\
$\Lg_5(\al)$ & $[e_1,e_2]=e_2, \,[e_1,e_3]=e_2+\al e_3,\,[e_1,e_4]=(\al +1)
e_4,\,[e_2,e_3]=e_4$\\
\end{tabular}
\end{center}   
\vspace*{0.5cm}
Here the decomposable algebras are not contained in the table.
They are given by: $\Lg_2(0)\cong
\Lr_{3,1}(\C)\oplus\C$, $\Lg_4(\al,0)\cong \Lr_{3,\al}(\C)\oplus\C$
with $\al \ne 0,1$ and $\Lg_4(0,1)\cong \Lr_3(\C)\oplus\C$. 
Note that $\Lg_5(\al)\cong \Lg_5(\al')$ if and only if $\al\al'=1$ or $\al=\al'$,
and $\Lg_4(\al,\be)\cong \Lg_4(\al',\be')$ if and only if the ratios
$1:\al:\be$ and $1:\al':\be'$ coincide (after some permutation).\\
We have given in \cite{BU10} the classification of all orbit closures in $\CL_4(\C)$.
One can improve the result as follows:  
\begin{prop}\label{1}
All degenerations in $\CL_4(\C)$ can be obtained by the composition
of the following essential degenerations:
$$
\Dg
\Lg_{4}(\alpha,\be) & \rTo & \Ln_4 & \rTo & \Ln_3\oplus \C & \rTo & \C^4 \\
\Lg_{4}(\alpha,1) & \rTo & \Lg_2(\alpha) & \rTo & \Ln_3\oplus \C  \\
\Lg_{4}(0,0) & \rTo & \Lr_2 \oplus \C & \rTo & \Ln_3\oplus \C \\
\Lg_2 (1) & \rTo & \Lg_1 & \rTo & \C^4 & \\
\Lg_5 (1) & \rTo & \Lg_3 & \rTo & \Lg_2 (2) & \\
\Lg_5 (\alpha) & \rTo & \Lg_4 (\alpha, \alpha +1) \\
\Ls\Ll_2\oplus\C & \rTo & \Lg_5(-1) \\
\Lr_2\oplus \Lr_2 & \rTo & \Lg_4(\alpha,0)\\
\Lr_2\oplus \Lr_2 & \rTo & \Lg_5(0) \\
\endDg
$$
\end{prop}
It is also possible to draw the diagrams of the degenerations in
the $4$ irreducible components. The diagram of the degenerations
in $\ov{O(\Lr_2\oplus \Lr_2)}$ looks as follows:

$$
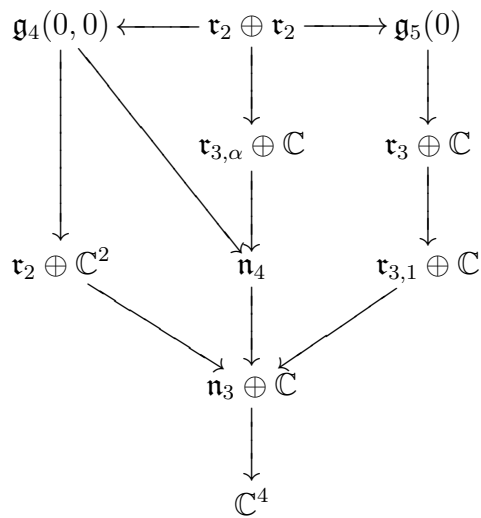
\Dg
\Lg_{4}(0,0) \aTo (2,-4) & \lTo & \Lr_2\oplus \Lr_2 & \rTo & \Lg_5(0) \\
\dTo         &      &  \dTo             &      &  \dTo \\
             &  \phantom{---}  & \Lr_{3,\alpha}\oplus\C & \phantom{--} & \Lr_3\oplus \C  \\
             &      & \dTo              &      &  \dTo \\
\Lr_2\oplus\C^2 \aTo (2,-2) &      &  \Ln_4            &      &  \Lr_{3,1}\oplus \C \aTo (-2,-2)\\
             &      &  \dTo             &      &  \\
             &      &  \Ln_3\oplus\C    &      &  \\
             &      &   \dTo            &      &  \\
             &      &   \C^4            &      &  \\
\endDg
$$              

To prove this classification result one uses the invariants mentioned
in proposition $\ref{inv}$. Moreover, of the algebra is solvable but
not nilpotent, the following numbers are of interest:
$$c_{ij}(\Lg)=\frac{\tr (\ad x)^i \tr (\ad y)^j}{\tr ((\ad x)^i\circ (\ad y)^j)}$$
If these numbers are independent of $x$ and $y$ in $\Lg$, and the denominator
does not vanish, then we obtain useful invariants.
For example, 
$$c_{ij}(\Lr_{3,\al}\oplus\C)=1+\frac{\al^i+\al^j}{1+\al^{i+j}}$$ 
In that case $c_{ij}(\Lh)=c_{ij}(\Lg)$ for all $\Lh\in \ov{O(\Lg)}$.

\end{document}